\documentclass[11pt]{article}
\usepackage{amsfonts}
\usepackage{amssymb}
\usepackage{amsthm}
\usepackage{eufrak}
\usepackage[dvips]{graphicx}
\usepackage{latexsym}
\newtheorem*{stat}{\hskip 0,65 cm \statname} %\unnumbered{stat}
\newenvironment{statement}[1]{\def\statname{#1}\medskip \begin{stat}}{\end{stat}}

\newtheorem{nstat}{\nstatname}[section]

\newtheorem{lemma}[nstat]{\hskip 0,65 cm Lemma}

\let\H\undefined
\let\L\undefined
\let\P\undefined
\newcommand{\Int}{\mathop{\mathrm{Int}}\nolimits}
\newcommand{\Bd}{\mathop{\mathrm{Bd}}\nolimits}
\newcommand{\id}{\mathop{\mathrm{id}}\nolimits}
\newcommand{\H}{\mathop{\cal{H}}\nolimits}
\newcommand{\M}{\mathop{\cal{M}}\nolimits}
\newcommand{\B}{\mathop{\cal{B}}\nolimits}
\newcommand{\P}{\mathop{\cal{P}}\nolimits}
\newcommand{\L}{\mathop{\cal{L}}\nolimits}
\newcommand{\I}{\mathop{\cal{I}}\nolimits}
\newcommand{\J}{\mathop{\cal{J}}\nolimits}
\let\phi\varphi

\begin{document}

\title{\large\bf LIFTING BRAIDS}
\author{\sc\normalsize M. Mulazzani\\
\sl\normalsize Dipartimento di Matematica\\[-3pt]
\sl\normalsize Universit\`a di Bologna -- Italy\\
\tt\small mulazza@dm.unibo.it \and
\sc\normalsize R. Piergallini\\
\sl\normalsize Dipartimento di Matematica e Fisica\\[-3pt]
\sl\normalsize Universit\`a di Camerino -- Italy\\
\tt\small pierg@camserv.unicam.it}
\date{}

\maketitle

\vskip 2 cm

\begin{abstract}

\medskip

\noindent In this paper we study the homeomorphisms of $B^2$ that
are liftable with respect to a simple branched covering. Since any
such homeomorphism maps the branch set of the covering onto itself
and liftability is invariant up to isotopy fixing the branch set,
we are dealing in fact with liftable braids.\\ We prove that the
group of liftable braids is finitely generated by liftable powers
of half-twists around arcs joining branch points. A set of such
generators is explicitly determined for the special case of
branched coverings $B^2 \to B^2$. As a preliminary result we also
obtain the classification of all the simple branched coverings of
$B^2$.

\medskip\smallskip\noindent
{\sl Key words and phrases\/}: branched covering of the disk,
liftable homeomorphism, liftable braid.

\smallskip\noindent
{\sl 2000 Mathematics Subject Classification\/}: 57M12.
\end{abstract}

\section*{Introduction}

\indent\indent A continuous map $p:F \to G$ between compact
surfaces with $G$ connected and oriented is a {\sl branched
covering} iff it is a local homeomorphism near $\Bd F$ and any
point $x \in \Int F$ has a neighborhood $U \subset \Int F$ such
that the restriction $p_{|U}:U \to p(U)$ is topologically
equivalent to the complex map $z \mapsto z^{d_x}$, for a uniquely
determined positive integer $d_x$,  the {\sl local order} of $p$
at $x$. In particular, we have $p(\Int F) = \Int G$ and $p(\Bd F)
= \Bd G$.

Given a branched covering $p:F \to G$, we denote by $S_p \subset
\Int F$ the (finite) set of the {\sl singular points} of $p$,
that is the points $x \in \Int F$ such that $d_x > 1$, and by
$B_p = \{P_1, \dots, P_n\} \subset \Int G$ the set of {\sl branch
points} of $p$, defined by $B_p = p(S_p)$. Then, the restriction
$p_|:F - p^{-1}(B_p) \to G - B_p$ is an ordinary covering with
$d$ sheets, where $d = d(p)$ is the {\sl order} of $p$. The
orientation of $G$ can be lifted to an orientation of $F$ which
makes $p:(F,\Bd F) \to (G, \Bd G)$ a map of positive degree
$d(p)$. {\sl We assume $F$ oriented in this way.}

Since $p_|$ uniquely determines $p$, by fixing a base point $*
\in G - B_p$ and numbering the fiber $p^{-1}(*)$, we can
represent $p$ by means of the monodromy $\phi_p: \pi_1(G - B_p,*)
\to \Sigma_d$ of the ordinary covering $p_|$, where $\Sigma_d$ is
the permutation group on $\{1,\dots,d\}$. We call $\phi_p$ the
{\sl monodromy} of $p$. In order to simplify the notation, we
write $\phi$ in place of $\phi_p$, when there is no risk of
confusion. Because of the choices of $\,*\,$ and of the numbering
of $p^{-1}(*)$, the monodromy is defined only up to inner
automorphisms of $\Sigma_d$.

A branched covering $p$ is called {\sl simple} iff it maps $S_p$
injectively onto $B_p$ and $d_x = 2$ for any $x \in S_p$. This
means that the monodromy of a small simple loop around any branch
point is a transposition.

\medskip

Two branched coverings $p:F \to G$ and $p':F' \to G'$ are {\sl
equivalent} iff there exist {\sl orientation preserving}
homeomorphisms $h:G \to G'$ and $k:F \to F'$ such that $p'k =
h\,p$.  Of course, in this case we have $d(p) = d(p')$, $h(B_p) =
B_{p'}$ and $k(S_p) = S_{p'}$. Now, it turns out that the
existence of a lifting $k:F \to F'$ of a given homeomorphism $h:G
\to G'$ such that $h(B_p) = B_{p'}$ depends only on the existence
of a lifting of the restriction $h_|:G-B_p \to G'-B_{p'}$. Then,
by the classical theory of ordinary covering, we get the following
criterion.

\medskip\begin{statement}{Lifting theorem} A homeomorphism
$h:G \to G'$ has a lifting  $k:F \to F'$ with respect to the
branched coverings $p:F \to G$ and $p':F' \to G'$ of the same
order $d$ iff $h(B_p) = B_{p'}$ and there exists an inner
automorphism $\varepsilon$ of $\Sigma_d$ such that $\phi_{p'} h_*
= \varepsilon\,\phi_p$, where $h_*:\pi_1(G - B_p,*) \to \pi_1(G' -
B_{p'},*')$ is the isomorphism induced by the restriction of $h$.
In this case $\varepsilon$ is given by the conjugation by $\sigma
= \nu'k\,\nu^{-1} \in \Sigma_d$, where $\nu: p^{-1}(*) \to
\{1,\dots,d\}$ and $\nu':p'^{-1}(*') \to \{1,\dots,d\}$ are the
numberings of the fibers $p^{-1}(*)$ and $p'^{-1}(*')$, with $*' =
h(*)$, inducing the monodromies $\phi_p$ and $\phi_{p'}$.
\end{statement}

\medskip

As an immediate consequence of this lifting theorem, we have an
equivalence criterion for branched coverings in terms of their
branch sets and monodromies.

\medskip\begin{statement}{Equivalence theorem}
Two branched coverings $p:F \to G$ and $p': F' \to G'$ of the same
order $d$ are equivalent iff there exist an orientation
preserving homeomorphism $h:G \to G'$ and an inner automorphism
$\varepsilon$ of $\Sigma_d$ such that $h(B_p) = B_{p'}$ and
$\phi_{p'}h_* = \varepsilon\,\phi_p$.
\end{statement}

\medskip

The classification of the simple branched coverings of $S^2$ up to
equivalence is classical and well known. In \cite{GK87} and
\cite{GK88}, Gabai and Kazez extended such classification to all
the closed  surfaces. The following Theorem A, giving a
classification of the simple branched coverings of $B^2$, is
stated without proof in \cite{BE84}. In Section
\ref{coverings/sec} we give a proof of Theorem A, by providing a
canonical way of representing branched coverings of $B^2$. We need
such canonical representation in order to get our main result
about liftable braids.

\medskip

Given a simple branched covering $p: F \to B^2$ of order $d$, we
fix a base point $* \in S^1$ and a numbering of $p^{-1}(*)$. Then,
we define {\sl total monodromy} of $p$ to be the permutation
$\phi_p(\omega)\in\Sigma_d$, where $\omega$ is the clockwise
oriented simple loop supported by $S^1$. Moreover, we denote by
$\Omega(p)$ the conjugation class of $\phi_p(\omega)$ in
$\Sigma_d$, which is uniquely determined by $p$ (actually by the
restriction of $p$ over $S^1$). Now we are in position to state
the classification theorem.

\medskip\begin{statement}{Theorem A}
Two connected simple branched coverings $p:F\to B^2$ and $p':
F'\to B^2$ are equivalent iff they have the same order $d$, the
same number $n$ of branch points and $\Omega(p) = \Omega(p')$.
\end{statement}

\medskip

Since $\Omega(p)$ is the class of $d$-cycles of $\Sigma_d$ for any
simple branched covering $p:F \to B^2$ with $\Bd F$ connected, by
the Riemann-Hurwicz formula we easily get the following corollary.

\medskip\begin{statement}{Corollary}
For every compact connected orientable surface $F$ with connected
boundary and for every integer $n \geq 2 - \chi(F)$ there exists
a unique (up to equivalence) simple covering  $p:F \to B^2$ with
$n$ branch points.
\end{statement}

\medskip

Given an orientable surface $S$ and a closed subset $C \subset S$,
we denote by $\H(S)$ the group of all the {\sl orientation
preserving} homeomorphisms of $S$ onto itself and by $\H(S,C)
\subset \H(S)$ the subgroup consisting of all the $h \in \H(S)$
such that $h(C) = C$. Moreover, if $D \subset S$ is another closed
subset, then we denote by $\H_D(S) \subset \H(S)$ and $\H_D(S,C)
\subset \H(S,C)$ the subgroups of the homeomorphisms which
coincide with the identity in $D$. Finally, we denote by $\M(S)$,
$\M(S,C)$, $\M_D(S)$ and $\M_D(S,C)$ the {\sl mapping class
groups} corresponding to the groups considered above (that is, we
set $\M = \pi_0\H$).

\medskip

For any $n \geq 1$, let $\B_n = \pi_1(\Gamma_n (\Int
B^2),\{P_1,\dots,P_n\})$ be the {\sl braid group} of order $n$ of
$S$ based at $\{P_1, \dots, P_n\} \subset \Int B^2$, where
$\Gamma_n(X) = (X^n - \Delta)/\Sigma_n$ denotes the configuration
space of all the subsets of $X$ with cardinality $n$.\break We
recall that there exists an isomorphism $\eta:\B_n \to
\M_{S^1}(B^2,\{P_1, \dots, P_n\})$, defined by setting $\eta(b)$
equal to the isotopy class of the ending homeomorphism $h_1$ of
any isotopy $t \mapsto h_t \in \H_{S_1}(B^2)$ which realizes the
braid $b$ (that is, the map $t \mapsto h_t(\{P_1, \dots, P_n\})$
is a loop in $\Gamma_n(\Int B^2)$ representing $b$).

\medskip

{\sl We use the right-handed notation for the action of braids on
everything, that is $(a)b$ denotes the image of $a$ by the action
of the braid $b$. If $a$ itself is a braid, then we have $(a)b =
b^{-1}ab$. Moreover, we adopt the following bracketing
convention:$(a)b_1b_2\dots b_n = (\dots((a)b_1)b_2\dots )b_n$.}

\medskip

We say that a homeomorphism $h \in \H_{\Bd G}(G)$ is {\sl
liftable} with respect to the branched covering $p:F \to G$ iff
there exists $k \in \H_{\Bd F}(F)$ such that $p\,k = h\,p$.\break
We call $k$ a {\sl lifting} of $h$. Of course, for such $h$ and
$k$, we have $h(B_p) = B_p$ and $k(S_p) = S_p$. Moreover, the
lifting $k$ is unique if $\Bd G \neq \emptyset$, otherwise $h$ may
have more than one lifting. In any case, liftability is preserved
by composition and is invariant by isotopy in $\H_{\Bd G}(G,B_p)$,
so it makes sense to speak of the (subgroup of the) liftable
isotopy classes in $\M_{\Bd_G}(G,B_p)$.

\medskip

Given a simple branched covering $p:F \to B^2$, we call $b \in
\B_n$ (the braid group based at the branch set $B_p = \{P_1,
\dots, P_n\}$ of $p$) a {\sl liftable braid} with respect to $p$
iff $\eta(b) \in \H_{S^1}(B^2,B_p)$ is a liftable isotopy class.
Moreover, we denote by $\L_p \subset \B_n$ the subgroup of the
liftable braids with respect to $p$.

\medskip

Following \cite{BW85}, we call {\sl curve} for the branched
covering $p:F \to B^2$ any simple arc $\alpha \subset B^2$ joining
the base point $* \in S^1$ with $B_p$ and such that $\Int\alpha
\subset \Int B^2 - B_p$. Curves are considered up to ambient
isotopy of $B^2$ modulo $S^1 \cup B_p$. A {\sl system of curves}
is any family of curves $\alpha_1, \dots, \alpha_k \subset B^2$
such that $\alpha_i \cap \alpha_j = \{*\}$ for all $i \neq j$. A
{\sl fundamental system} is a maximal system of curves, that is a
system of curves $\alpha_1, \dots, \alpha_n$ with the same
cardinality of $B_p = \{P_1, \dots, P_n\}$.

For any curve $\alpha$, let $\lambda_\alpha \in \pi_1(B^2 -
B_p,*)$ be the homotopy class of a simple loop supported by the
clockwise oriented boundary of a small regular neighborhood
$N(\alpha)$ of $\alpha$ in $B^2$. In order to simplify the
notation we put $\phi_p(\alpha) = \phi_p(\lambda_\alpha)$. If $p$
is simple, then $\phi_p(\alpha)$ is a transposition for any curve
$\alpha$. Viceversa, if $\alpha_1, \dots, \alpha_n$ is a
fundamental system for $p$ and $\phi_p(\alpha_i)$ is a
transposition for any $i = 1, \dots, n$, then $p$ is simple.

Since $\lambda_{\alpha_1}, \dots, \lambda_{\alpha_n}$ generate
$\pi_1(B^2 - B_p,*)$, for any fundamental system $\alpha_1,
\dots,\\ \alpha_n$ for $p$, the branched covering $p$ is
completely determined, up to equivalence, by the monodromies
$\phi_p(\alpha_1), \dots, \phi_p(\alpha_n)$.

\medskip

Following \cite{BW85} again, we call {\sl interval} for the
branched covering $p:F \to B^2$ any simple arc $x \subset B^2$
joining two branch points and such that $\Int x \subset \Int B^2 -
B_p$. Intervals are considered up to ambient isotopy of $B^2$
modulo $S^1 \cup B_p$.\break We call interval and we denote by the
same symbol $x$ also the counterclockwise half-twist around $x$
and the corresponding braid in $\B_n$.

It immediately follows from the lifting theorem that any interval
$x$ has a liftable power. In fact, we prove in Section
\ref{braids/sec} that, if $x$ is not liftable, then either $x^2$
or $x^3$ is liftable.

\medskip\begin{statement}{Theorem B}
For any branched covering $p:F \to B^2$, the group of the liftable
braids $\L_p$ is finitely  generated by liftable powers of
intervals.
\end{statement}

\medskip

The proof of Theorem B is given in Section \ref{generalcase/sec}.
As a preliminary step, in Section \ref{specialcase/sec} we
consider the special case of $F = B^2$. In this case, we
explicitly provide a set of generators, as described in the
following Theorem C.

\medskip

Let $p_n: B^2 \to B^2$ be the unique (up to equivalence) simple
covering of order $d=n+1$ with $n$ branch points. For sake of
simplicity, we denote by $\L_n \subset \B_n$ the group $\L_{p_n}$
of the liftable braids respect to $p_n$.

We assume the fundamental system $\alpha_1,\dots, \alpha_n$ and
the numbering of the sheets of $p_n$ be fixed, in such a way that
the sequence of transpositions $\phi(\alpha_1), \dots,\\
\phi(\alpha_n)$ is in the {\sl canonical form} $(1\;2), \dots,
(d{-}1\;d)$, that is $\phi(\alpha_i) = (i\;i{+}1)$ for each $i =
1, \dots, n$.

For each $i=1, \dots, n-1$, we define $x_i \simeq \alpha_i \cup
\alpha_{i+1}$ to be the unique interval such that $x_i \cup
\alpha_i \cup \alpha_{i+1}$ is a Jordan curve whose interior does
not contain any branch point. Moreover, we put $x_{i,j} =
(x_i)x_{i+1} \dots x_{j-1}$, for $1 \leq i < j \leq n$.

\medskip\begin{statement}{Theorem C}
For any $n > 1$, the group $\L_n$ of the liftable braids with
respect to the branched covering $p_n:B^2 \to B^2$ is generated
by the above described braids $x_i^3$ and $x_{i,j}^2$, with $1
\leq i < n$ and $i + 1 < j \leq n$.
\end{statement}

\medskip

The above theorems constitute the first results in the study of
the lifting homomorphism $\lambda_p:\L_p \to \M_{\Bd F}(F)$, we
are planning to carry out in order to find a set of normal
generators of $\ker\lambda_p$ in $\L_p$, for any branched covering
$p:F \to B^2$. This would generalize a result obtained in
\cite{BW85} (see also \cite{BW94}) for coverings of degree 3.

Our work is mainly aimed to get a general equivalence theorem for
simple branched coverings of $S^3$ in terms of covering moves. In
fact, the equivalence theorem for degree 3 coverings given in
\cite{P91} and \cite{P95} is essentially based on the results of
\cite{BW85}.

\section{Branched coverings of $B^2$\label{coverings/sec}}

\indent\indent Let $p: F \to B^2$ be a simple branched covering of
order $d$. Given a fundamental system $\alpha_1, \dots, \alpha_n$
with monodromies $\phi(\alpha_i) = (j_i\;k_i)$ for $i = 1, \dots,
n$, we define the non-oriented graph $\Gamma = \Gamma_p(\alpha_1,
\dots, \alpha_n)$ to have vertices $v_1, \dots, v_d$ and edges
$e_1 = \{v_{j_1}, v_{k_1}\}, \dots, e_n=\{v_{j_n}, v_{k_n}\}$.

We remark that the ordering of the vertices of $\Gamma$ is not
relevant, since it depends on an arbitrary numbering of the
sheets of $p$. On the contrary, the ordering of the edges
contains relevant information related to the choice of the
fundamental system. Therefore, we consider $\Gamma$ as an {\sl
edge-ordered graph}, in such a way that $p$ is uniquely
determined by $\Gamma$ up to equivalence.

Moreover, for each non-oriented edge-ordered graph $\Gamma$ with
$d$ vertices and $n$ edges, there exist a simple branched
covering $p = p_\Gamma : F_\Gamma \to B^2$ of order $d$ with $n$
branch points and a fundamental system $\alpha_1, \dots, \alpha_n$
for $p$, such that $\Gamma = \Gamma_p(\alpha_1, \dots, \alpha_n)$.

\medskip\begin{lemma}\label{homotopy/thm} % Lemma 1.2.
$F_\Gamma$ has the same homotopy type of\/ $\Gamma$.
\end{lemma}

\medskip\begin{proof} $F_\Gamma$ is homeomorphic to the topological union $D_1 \sqcup
\dots \sqcup D_d$ of $d$ discs, with a band glued between
$D_{j_i}$ and $D_{k_i}$ for every $i = 1,\dots,n$.
\end{proof}

Now we want to establish when two connected edge-ordered graphs
$\Gamma$ and $\Gamma'$, with $d$ vertices (which can be assumed
to be the same) and $n$ edges, determine equivalent coverings
$p_\Gamma$ and $p_{\Gamma'}$.

For any $i=1, \dots, n$, let $O_i$ be the {\sl elementary move}
corresponding to the transformation (3.6) of \cite{BE79}, which
transforms the graph $\Gamma$ with edges $e_1, \dots, e_n$ in the
graph $\Gamma'$ with edges $e'_1, \dots, e'_n$, defined in the
following way: if $e_i$ and $e_{i+1}$ are disjoint or they share
both the endpoints, then we put $e'_i = e_{i+1}$, $e'_{i+1} =
e_i$ and $e'_k = e_k$ for $k \neq i,i+1$; otherwise, if $e_i$ and
$e_{i+1}$ share only one endpoint, say $e_i = \{a,b\}$ and
$e_{i+1} = \{b,c\}$ with $a \neq c$, then we put $e'_i = \{a,c\}$,
$e'_{i+1} = e_i$ and $e'_k = e_k$ for $k \neq i,i+1$. We also
denote by $O_i^{-1}$ the inverse elementary moves, defined in the
obvious way.

Fixed a numbering of the vertices $v_1, \dots, v_d$ of $\Gamma$,
we associate to each edge $e_i = \{v_{j_i},v_{k_i}\}$ the
transposition $\tau_i = (j_i\;k_i) \in \Sigma_d$. Then, we define
$\Omega(\Gamma)$ as the conjugation class of the product $\tau_1
\cdots \tau_n$ in $\Sigma_d$. We observe that the class
$\Omega(\Gamma)$ is uniquely determined by $\Gamma$ (that is it
does not depend on the numbering of the vertices) and is
invariant with respect to elementary moves. Furthermore, it is
straightforward to see that $\Omega(\Gamma) = \Omega(p_\Gamma)$.

\medskip\begin{lemma} \label{graph/thm} % Lemma 1.3.
Let\/ $\Gamma$ and\/ $\Gamma'$ be two connected edge-ordered
graphs with $d$ vertices and $n$ edges. Then the coverings
$p_\Gamma$ and $p_{\Gamma'}$ are equivalent if and only if\/
$\Omega(\Gamma) = \Omega(\Gamma')$.
\end{lemma}

\medskip\begin{proof}
The `only if' part is trivial. Viceversa, it suffices to prove
that each connected edge-ordered graph $\Gamma$ with $d$ vertices
and $n$ edges can be transformed, by using elementary moves and
their inverses, into a canonical form dependent only on
$\Omega(\Gamma)$.

Let us denote by $c_1 \geq \cdots \geq c_m$ the cardinalities of
the non-trivial orbits generated by any permutation of
$\Omega(\Gamma)$ and let $l_i = c_1 + \cdots + c_i$ for each $i =
1, \dots, m$, then we can choice as a canonical representative of
$\Omega(\Gamma)$ the permutation $\pi$ given by the product
$$(1\;2) \cdots (l_1{-}1\;l_1)(l_1{+}1\;l_1{+}2) \cdots
(l_2{-}1\;l_2) \ \cdots \
         (l_{m-1}{+}1\;l_{m-1}{+}2) \cdots (l_m{-}1\;l_m)\,.$$

On the other hand, there exists a numbering $v_1, \dots, v_d$ of
the vertices of $\Gamma$, such that $\tau_1 \cdots \tau_n = \pi$,
where $\tau_i = (j_i\;k_i)$ is the transposition associated to
the edge $e_i = \{v_{j_i}, v_{k_i}\}$, for every $i = 1, \dots,
n$.

We want to transform $\Gamma$ by elementary moves, leaving the
numbering of the vertices fixed, in such a way that the sequence
$\tau_1, \dots, \tau_n$ becomes
$$\matrix{
(1\;2), \dots, (l_1{-}1\;l_1),(l_1\;l_1{+}1),(l_1\;l_1{+}1),
   (l_1{+}1\;l_1{+}2), \dots, (l_2{-}1\;l_2),(l_2\;l_2{+}1),\cr
(l_2\;l_2{+}1),\dots,(l_{m-1}\;l_{m-1}{+}1),(l_{m-1}\;l_{m-1}{+}1),
   (l_{m-1}{+}1\;l_{m-1}{+}2), \dots, (l_m{-}1\;l_m),\cr
(l_m\;l_m{+}1),(l_m\;l_m{+}1),(l_m{+}1\;l_m{+}2),(l_m{+}1\;l_m{+}2),
\dots,
   (d{-}1\;d),(d{-}1\;d),\cr
(d{-}1\;d),(d{-}1\;d), \dots,
(d{-}1\;d),(d{-}1\;d)\;,\kern-1ex\cr}$$ where the first two rows
contain the transposition sequence defining $\pi$, with
additional pairs of consecutive equal transpositions inserted
between disjoint cycles, and the last two rows consist of pairs
of equal consecutive transpositions. Moreover: a) if $\pi$ is the
identity then the first two rows are empty and $l_m = m = 1$; b)
if $l_m = d$ then the third row is empty; c) the fourth row
contains $(n - m + l_m)/2 - d + 1$ pairs of transpositions.

We proceed by induction on $n$. If $n = 1$, then $\Gamma$ itself
has the required form, in fact the only possibility is $\tau_1 =
(1\;2)$ and $d = 2$, since $\Gamma$ is connected. In the rest of
the proof, we deal with the inductive step, assuming $n > 1$.

To begin with, we show how to perform moves on $\Gamma$ in order
to obtain a sequence $\tau_1, \dots, \tau_n$ of the type
$(j_1\;k_1), \dots, (j_{n'}\;k_{n'}), (d{-}1\;d), \dots,
(d{-}1\;d)$, with $0 \leq n' < n$ and $j_i,k_i \neq d$ for each
$i \leq n'$. First of all, by using Remark (3.7) of \cite{BE79},
it is easy to get a sequence having the form $(j_1\;k_1), \dots,
(j_{n'}\;k_{n'}), (j_{n'+1}\;d), \dots, (j_n\;d)$, with $0 \leq
n' < n$ and $j_i,k_i \neq d$ for each $i \leq n'$. Then, since
$O_i$ change the pair $(j_i\;d),(j_{i+1}\;d)$ with $j_i \neq
j_{i+1}$ into the pair $(j_i\;j_{i+1}), (j_i\;d)$, we can limit
ourselves to consider only the case $j_{n'+1} = \dots = j_n =
h$.  If $h = d-1$, we have done. Otherwise, if $h \neq d-1$, we
have that $n' > 0$, by connectedness, and that $n - n'$ is even,
since $(d)\tau_1 \dots \tau_n$ can only assume the value $d-1$ or
$d$. At this point, we could get the desired form by the sequence
of elementary moves $O_{n'},\dots, O_{n-1}, O_{n-1},\dots,
O_{n'}$ if $(j_{n'}\;k_{n'})$ was $(h\;d{-}1)$. So, it remains to
show how to obtain $e_{n'} = \{v_h,v_{d-1}\}$ without changing the
edges $e_i$ with $i > n'$. By connectedness, there exists a chain
of edges $e_{i_1}, \dots, e_{i_l}$ of minimum length $l " 1$
connecting $v_h$ and $v_{d-1}$, with $i_1, \dots, i_l \leq n'$.
If $l = 1$, then $\tau_i = (h\;d{-}1)$ and we can finish by using
Remark (3.7) of \cite{BE79} again. If $l > 1$ and $\vert i_l -
i_{l-1} \vert = 1$, then we can decrease by one the length of the
chain by performing the move $O_i$ with $i = \min(i_l, i_{l-1})$.
On the other hand, if $\vert i_l - i_{l-1} \vert > 1$, then we
can reduce by one the difference between $i_{l-1}$ and $i_l$, by
performing either $O_{i_l - 1}^{-1}$ if $i_{l-1} < i_l$ or
$O_{i_l + 1}$ if $i_{l-1} > i_l$. So we can conclude this part of
the proof by induction on $l$ and $\vert i_l - i_{l-1} \vert$.

 From now on, we assume that the first $n'$ edges of $\Gamma$ do not
contain $v_d$ and that all the last $n - n'$ edges of $\Gamma$
join $v_{d-1}$ and $v_d$. If $n' = 0$, then we have finished ($d
= 2$ and either $\pi = \id$ or $\pi=(1\;2)$ depending on the
parity of $n$).

Let us consider the case $n' > 0$. We denote by $\Gamma' \subset
\Gamma$ the subgraph having vertices $v_1, \dots, v_{d-1}$ and
edges $e_1, \dots, e_{n'}$. Since $\Gamma'$ is connected and the
permutation $\tau_1 \dots \tau_{n'}$ is of the type requested for
$\pi$, we can apply the induction hypothesis, in order to
transform $\Gamma'$ into the canonical form, by a sequence of
elementary moves and inverse of them. The same sequence of moves
also transforms $\Gamma$ in a canonical form, possibly except for
the presence of more than two transpositions $(d{-}2\;d{-}1)$
immediately before the $n-n'$ transpositions $(d{-}1\;d)$. In
fact, the canonical form for $\Gamma$ contains either one
transposition $(d{-}2\;d{-}1)$ if $c_m > 2$ or two of them if
$c_m = 2$. Hence, to complete the proof, it suffices to change
all the $n''$ exceeding transpositions $(d{-}2\;d{-}1)$ into
$(d{-}1\;d)$. Taking into account that such transpositions are
preceded by at least one more $(d{-}2\;d{-}1)$ and followed by
$(d{-}1\;d)$ and that their number $n''$ is even, we can realize
the wanted change by the sequence of elementary moves $O_{n'},
\dots, O_{n'-n''+1}, O_{n'-n''+1}, \dots, O_{n'}, O_{n'-n''},
\dots, O_{n'-1}, O_{n'-1}, \dots, O_{n'-n''}$.
\end{proof}

At this point, we can prove the Theorem A stated in the
introduction.

\medskip \begin{proof}[Proof of Theorem A]
Let $\Gamma$ and $\Gamma'$ two edge-ordered graphs such that
$p=p_\Gamma$ and $p'=p_{\Gamma'}$. By Lemma \ref{homotopy/thm},
$\Gamma$ and $\Gamma'$ are both connected. Moreover, we have
$\Omega(\Gamma) = \Omega(p) = \Omega(p') = \Omega(\Gamma')$.
Therefore, Theorem A immediately follows from Lemma
\ref{graph/thm}.
\end{proof}

We conclude this section by considering some elementary
properties of the restrictions of a covering of $B^2$ over
subdisks, which we will need in the following sections.

\medskip

Given a simple branched covering $p:F \to B^2$ with $n$ branch
points and a system of curves $\alpha_1, \dots, \alpha_k \subset
B^2$ for $p$, we denote by $p^{\alpha_1, \dots, \alpha_k}:
F^{\alpha_1, \dots, \alpha_k} \to B^2_{\alpha_1, \dots, \alpha_k}$
the restriction of $p$ to $F^{\alpha_1, \dots, \alpha_k} =
p^{-1}(B^2_{\alpha_1, \dots, \alpha_k})$, where $B^2_{\alpha_1,
\dots, \alpha_k}$ is the disk $B^2 - \Int N(\alpha_1, \dots,
\alpha_k)$, being $N(\alpha_1, \dots, \alpha_k)$ a regular
neighborhood of $\alpha_1 \cup \dots \cup \alpha_k$ that does not
contain any branch points other than the endpoints of the curves
$\alpha_1, \dots, \alpha_k$.

We remark that $p^{\alpha_1, \dots, \alpha_k}$ is a simple
covering uniquely determined up to equivalence, which has the
same order of $p$ and $n-k$ branch points. Moreover, if $p$ has
$c$ components, $p^{\alpha_1, \dots, \alpha_k}$ has $c'$
components, with $c\le c'\le c+k$.

As base point for $B^2_{\alpha_1, \dots, \alpha_k}$ we can choose
either the starting-point $*'$ or the ending-point $*''$ of the
arc $N(\alpha_1, \dots, \alpha_k) \cap S^1$, oriented accordingly
with the counterclockwise orientation of $S^1$. We denote by
$\omega'_{\alpha_1, \dots, \alpha_k}$ (resp. $\omega''_{\alpha_1,
\dots, \alpha_k}$) the simple clockwise oriented loop based at
$*'$ (resp. $*''$) and supported by $\Bd B^2_{\alpha_1, \dots,
\alpha_k}$. The liftings of the arc $N(\alpha_1, \dots, \alpha_k)
\cap S^1$ determine bijections $p^{-1}(*) \cong p^{-1}(*') \cong
p^{-1}(*'')$. By means of this bijections, any numbering of the
sheets of $p$ induces a {\sl coherent numbering} of the sheets of
$p^{\alpha_1, \dots, \alpha_k}$, depending on the choice of $*'$
or $*''$ as base point $B^2_{\alpha_1, \dots, \alpha_k}$.

\medskip

{\sl The sheets of any restriction $p^{\alpha_1, \dots,
\alpha_k}$ of $p$ will be ever numbered coherently with the ones
of $p$, whichever will be the choice of the base point for
$B^2_{\alpha_1, \dots, \alpha_k}$. We will denote with the same
letter $\phi$ the monodromy of $p^{\alpha_1, \dots, \alpha_k}$,
with respect to this coherent numbering, as well as the monodromy
of $p$, without any explicit reference to the choice of the base
point.}

\medskip

Given any curve $\alpha \subset B^2$ for $p$ such that $\alpha
\cap (\alpha_1 \cup \dots \cup \alpha_k) = \{*\}$, we can assume
(up to isotopy) that $\alpha \cap B^2_{\alpha_1, \dots,
\alpha_k}$ is an arc. Then, we denote by $\alpha'$ and (resp.
$\alpha''$) the curve in $B^2_{\alpha_1, \dots, \alpha_k}$
obtained by sliding the initial point of $\alpha \cap
B^2_{\alpha_1, \dots, \alpha_k}$ along the arc $N(\alpha_1,
\dots, \alpha_k) \cap B^2_{\alpha_1, \dots, \alpha_k}$ until $*'$
(resp. $*''$) is reached. By using these notations we can write
$p^{\alpha_1, \dots, \alpha_k} \cong (p^{\alpha_1, \dots,
\alpha_h})^{\alpha_{h+1}'', \dots, \alpha_k''} \cong
(p^{\alpha_{h+1}, \dots, \alpha_k})^{\alpha_1', \dots, \alpha_h'}$
for each $h = 1, \dots, k-1$.

\medskip

If $\alpha_1,\dots,\alpha_n$ is a fundamental system for $p$, then
for any $i_1 < \dots < i_k$ and $j_1 < \dots < j_{n-k}$ such that
$\{i_1, \dots, i_k\} \cup \{j_1, \dots, j_{n-k}\}=\{1,\dots,n\}$,
we have that $\alpha_{j_1}',\dots,\alpha_{j_{n-k}}'$ (resp.
$\alpha_{j_1}'', \dots, \alpha_{j_{n-k}}''$) is a fundamental
system for $p^{\alpha_{i_1}, \dots, \alpha_{i_k}}$ with base
point $*'$ (resp. $*''$). Moreover, by putting $\tau_i =
\phi(\alpha_i)$, straightforward computations give the following
equalities: $\phi(\alpha_j') = \tau_j$ and $\phi(\alpha_j'') =
(\tau_j)\tau_{i_1} \dots \tau_{i_k}$, if $j < i_1$;
$\phi(\alpha_j') = (\tau_j)\tau_{i_l}\dots \tau_{i_1}$ and
$\phi(\alpha_j'') = (\tau_j)\tau_{i_{l+1}} \dots \tau_{i_k}$, for
$i_l < j < i_{l+1}$ with $1 \leq l \leq k-1$; $\phi(\alpha_j') =
(\tau_j)\tau_{i_k}\dots \tau_{i_1}$ and $\phi(\alpha_j'') =
\tau_j$, if $j > i_k$.

\medskip \begin{lemma} \label{restriction/thm} % Lemma 1.1.
If $p:F \to B^2$ is a simple covering and $\alpha_1, \dots,
\alpha_k \subset B^2$ is a system of curves for $p$, then
$\phi(\omega'_{\alpha_1,\dots,\alpha_k})= \phi(\omega)
\phi(\alpha_k) \dots \phi(\alpha_1)$ and $\phi(\omega''_{\alpha_1,
\dots, \alpha_k}) = \phi(\alpha_k) \dots
\phi(\alpha_1)\phi(\omega)$.
\end{lemma}

\medskip\begin{proof}
Let $\alpha_{k+1}, \dots, \alpha_n$ be curves such that $\alpha_1,
\dots, \alpha_n$ is a fundamental system for $p$. Then, by the
equalities above, $\phi(\omega'_{\alpha_1, \dots, \alpha_k}) =
\phi(\alpha_{k+1}') \dots \phi(\alpha_n') = (\phi(\alpha_{k+1})
\dots \phi(\alpha_n)) \phi(\alpha_k) \dots \phi(\alpha_1) =
\phi(\omega) \phi(\alpha_k) \dots \phi(\alpha_1)$ and
$\phi(\omega_{\alpha_1,\dots,\alpha_k}'') =\\ \phi(\alpha_{k+1}'')
\dots \phi(\alpha_n'') = \phi(\alpha_{k+1}) \dots \phi(\alpha_n) =
\phi(\alpha_k) \dots \phi(\alpha_1)\phi(\omega)$.
\end{proof}

\medskip \begin{lemma} \label{disk/thm} % Lemma 1.7.
A connected simple covering $p:F \to B^2$ with $n$ branch points
is equivalent to $p_n$ if and only if one of the following
conditions holds: $(1)$ $F \cong B^2$; $(2)$ $p$ has order $n+1$;
$(3)$ $p^\alpha$ is disconnected for every curve $\alpha$.
\end{lemma}

\medskip\begin{proof}
$p \cong p_n \Rightarrow (1)$ is trivial. $(1) \Rightarrow (2)$
follows from Lemma \ref{homotopy/thm}.\break $(2) \Rightarrow
(3)$ follows from the fact that $n-1$ transpositions cannot
generate a transitive action on $\{1, \dots, n+1\}$. In order to
prove the implication $(3) \Rightarrow p \cong p_n$, we consider
a fundamental system $\alpha_1, \dots, \alpha_n$ for $p$ such
that the sequence of transpositions $\phi_p(\alpha_1), \dots,
\phi_p(\alpha_1)$ has the canonical form obtained in the proof of
Lemma \ref{graph/thm}. It is easy to see that $p \cong p_n$ iff
$\phi(\alpha_m) \neq \phi(\alpha_{m+1})$ for each $m = 1, \dots,
n-1$. On the other hand, if $\phi(\alpha_m) = \phi(\alpha_{m+1})$
for some $m$, then the restriction $p^{\alpha_m}$ is connected.
\end{proof}

\section{Liftable braids and intervals\label{braids/sec}}

\indent\indent In this section we consider a simple branched
covering $p:F \to B^2$ of degree $d$ with $n$ branch points and
denote by $\phi$ its monodromy. We denote by $\B_n$ the braid
group based at the branch set $B_p$ of $p$ and by $\L_p \subset
\B_n$ the subgroups of the liftable braids with respect to $p$.

Let us begin with some elementary properties of liftable braids.
The following liftability criterion in terms of action on a
fundamental system will play a crucial role.

\medskip \begin{lemma} \label{liftability/thm} % Lemma 1.8.
Let $\alpha_1, \dots, \alpha_n$  be a fundamental system for $p$.
Then, a braid $b \in \B_n$ is liftable if and only if
$\phi((\alpha_i)b) = \phi(\alpha_i)$ for every $i = 1, \dots, n$.
\end{lemma}

\medskip\begin{proof}
Let $b_*$ be the automorphism of $\pi_1(B^2 - B_p,*)$ induced by
the restriction of $b$ to $B^2- B_p$. By the lifting theorem, $b$
is liftable if and only if $\phi\,b_* = \phi$, since the lifting
of $b$ is the identity on $p^{-1}(*)$ and thus it induces the
identity conjugation on $\Sigma_d$. Then, the statement follows
from the fact that the sequence of transpositions associated to a
fundamental system uniquely determines $\phi$.
\end{proof}

Let $\alpha_1, \dots, \alpha_k$ and $\beta_1, \dots, \beta_k$ two
systems of curves. Suppose that there exists a liftable braid
$b\in \B_n$ such that $(\alpha_i)b = \beta_i$ for every $i =
1,\dots,k$. Since any system of curves can be completed to a
fundamental system, by the previous lemma, we have $\phi(\alpha_i)
= \phi(\beta_i)$ for every $i = 1, \dots, k$. On the other hand,
we can always suppose that $(B^2_{\alpha_1, \dots, \alpha_k})b =
B^2_{\beta_1, \dots, \beta_k}$, up to isotopy. Therefore, the
restriction $b_|: B^2_{\alpha_1, \dots, \alpha_k} \to
B^2_{\beta_1, \dots, \beta_k}$ induces an equivalence between
$p^{\alpha_1, \dots, \alpha_k}$ and $p^{\beta_1, \dots,
\beta_k}$, preserving the numbering of the sheets, when they are
referred to the same base point $*'$ or $*''$.

\medskip \begin{lemma} \label{systems/thm} % Lemma 1.9.
Given two systems of curves $\alpha_1, \dots, \alpha_k$ and
$\beta_1, \dots, \beta_k$ for $p$, there exists a liftable braid
$b\in \B_n$ such that $(\alpha_i)b = \beta_i$ for every $i = 1,
\dots, k$ if and only if the following conditions hold: $(1)$
$\phi(\alpha_i) = \phi(\beta_i)$ for every $i = 1,\dots, k$; $(2)$
there exists a bijection between the components of
$p^{\alpha_1,\dots, \alpha_k}$ and $p^{\beta_1, \dots, \beta_k}$
preserving the number of branch points and the numbering of the
sheets, when they are referred to the same base point.
\end{lemma}

\medskip\begin{proof}
The `only if' part immediately follows from the previous
discussion. In order to prove the converse, it suffices to extend
the given systems of curves to fundamental systems $\alpha_1,
\dots, \alpha_n$ and $\beta_1, \dots, \beta_n$ such that
$\phi(\alpha_i) = \phi(\beta_i)$ for every $i = 1, \dots, n$. In
fact, in this case the braid $b \in B_n$ uniquely defined by
$(\alpha_i)b = \beta_i$ for $i = 1, \dots, n$ turns out to be
liftable by Lemma \ref{liftability/thm}. The fundamental systems
$\alpha_1, \dots, \alpha_n$ and $\beta_1, \dots, \beta_n$ will be
constructed by induction on $m = n-k$.

If $m = 0$ there is nothing to do. So, we assume $m > 0$ and
observe that in this case there exist connected components $C
\subset F^{\alpha_1, \dots, \alpha_k}$ and $D \subset F^{\beta_1,
\dots, \beta_k}$, such that the restrictions $p^{\alpha_1, \dots,
\alpha_k}_{|C}:C \to B^2_{\alpha_1, \dots, \alpha_k}$ and
$p^{\beta_1, \dots, \beta_k}_{|D}:D \to B^2_{\beta_1, \dots,
\beta_k}$ are non-trivial. We assume that such components
correspond each other with respect to the bijection of property
(2). Then, they involve the same sheets $\{i_1, \dots, i_e\}$ and
the same number $l > 0$ of branch points. Moreover, by Lemma
\ref{restriction/thm}, they have the same total monodromy with
respect to the base point $*''$, that is the restriction of
$\phi(\omega'' _{\alpha_1, \dots, \alpha_k}) = \phi(\omega''
_{\beta_1, \dots, \beta_k})$ to $\{i_1, \dots, i_e\}$.

By the proof of Theorem A, it is possible to construct two
fundamental systems $\gamma_1, \dots, \gamma_l$ for $p^{\alpha_1,
\dots, \alpha_k}_{|C}$ and $\delta_1, \dots, \delta_l$ for
$p^{\beta_1, \dots, \beta_k}_{|D}$, with the same base point
$*''$ and such that $\phi(\gamma_i) = \phi(\delta_i)$ for every
$i = 1, \dots, l$.

Now we consider the systems of curves $\alpha_1, \dots,
\alpha_{k+l}$ and $\beta_1, \dots, \beta_{k+l}$, extending the
original ones in such a way that $\alpha''_i = \gamma_{i-k}$ and
$\beta''_i = \delta_{i-k}$ for all $i = k+1, \dots, k+l$.
Properties (1) and (2) still hold for these new systems of
curves. Therefore, by the induction hypothesis, they can be
further extended to fundamental systems as desired.
\end{proof}

We remark that property (2) in the statement of Lemma
\ref{systems/thm} trivially follows from property (1) when the
restrictions $p^{\alpha_1, \dots, \alpha_k}$ and $p^{\beta_1,
\dots, \beta_k}$ are connected.  More generally, this fact holds
also when the two restrictions have at most one non-trivial
component and their trivial sheets are numbered in the same way.

\medskip

In the rest of this section, we deal with intervals. Given an
interval $x \subset B^2$ for the covering $p$, we say that $x$ is
of {\sl type $i$} iff $x^i$ is the first positive power of $x$
which is liftable with respect to $p$ as a braid.

By the following lemma (cf. Lemma 2.4 of \cite{BW85}), each
interval is either of type 1 or type 2 or type 3. Moreover, it
can be easily realized that the intervals $x$ and $(x)b$ are of
the same type for each liftable braid $b \in \L_p$.

\medskip \begin{lemma} \label{intervals/thm} %Lemma 1.10.
Let $x$ be an interval for $p$ and $\alpha$ be a curve for $p$
meeting $x$ only at one of its endpoints. Then: $x$ is of type
$1$ if and only if $\phi(\alpha) = \phi((\alpha)x)$; $x$ is of
type $2$ if and only if $\phi(\alpha)$ and $\phi((\alpha)x)$ are
disjoint transpositions; $x$ is of type $3$ if and only if
$\phi(\alpha)$ and $\phi((\alpha)x)$ are different and not
disjoint.
\end{lemma}

\medskip\begin{proof}
Given $x$ and $\alpha$ as in the statement, let $\alpha_1, \dots,
\alpha_n$ be a fundamental system such that $\alpha_1 = \alpha$,
$\alpha_2 = (\alpha)x$ and $(\alpha_i)x = \alpha_i$ for $i = 3,
\dots, n$. By Lemma \ref{liftability/thm}, $x$ is liftable iff it
preserves all the monodromies of such fundamental system, that is
iff $\phi(\alpha_1) = \phi(\alpha_2)$. The other two cases can be
achieved by similar applications of Lemma \ref{liftability/thm}
to the intervals $x^2$ and to $x^3$, taking into account that
$\phi((\alpha)x^2) = \phi((\alpha)x) \phi(\alpha) \phi((\alpha)x)$
and $\phi((\alpha)x^3) = \phi((\alpha)x^2) \phi((\alpha)x)
\phi((\alpha)x^2) = \phi((\alpha)x) \phi(\alpha) \phi((\alpha)x)
\phi(\alpha) \phi((\alpha)x)$.
\end{proof}

We denote by $\I_p \subset \L_p$ the subgroup generated by all the
liftable powers of intervals, that is by the intervals of type 1,
the second power of the intervals of type 2 and the third power
of the intervals of type 3. Of course, Theorem B says that $\I_p
= \L_p$. Nevertheless, it is temporarily convenient to keep
different notations for the two groups.

\medskip

Fixed a fundamental system $\alpha_1, \dots, \alpha_n$ for $p$, we
call {\sl index} of a curve or an interval (with respect to
$\alpha_1, \dots, \alpha_n$) the minimum number (up to isotopy)
of the intersections with $\alpha_1 \cup \dots \cup \alpha_n$, not
including the endpoints.

Moreover, depending on the fundamental system $\alpha_1, \dots,
\alpha_n$, we give the following definitions: $x_i \simeq
\alpha_i \cup \alpha_{i+1}$ is the unique interval such that $x_i
\cup \alpha_i \cup \alpha_{i+1}$ is a Jordan curve whose interior
does not contain any branch point, for $i = 1, \dots, n-1$;
$x_{i,j} = (x_i)x_{i+1} \dots x_{j-1}$, for $1 \leq i < j \leq
n$; $\widehat x_{i,j} = (x_i)x_{i+1}^{-1} \dots x_{j-1}^{-1}$,
for $1 \leq i < j \leq n$; in addition, as a notational
convenience, we put $x_{i,j} = x_{j,i}$ and $\widehat x_{i,j} =
\widehat x_{j,i}$, for $1 \leq j < i \leq n$. In particular, we
have $x_i = x_{i,i+1} = \widehat x_{i,i+1}$.

We remark that the braids $x_1,\dots,x_{n-1}$ are the usual
standard generators of the braid group $\B_n$; similarly, the
braids $x_{i,j}^2$ (as well as the braids $\widehat x_{i,j}^2$)
with $1 \leq i < j \leq n$ are standard generators of the pure
braids group $\P_n \subset \B_n$.

\medskip

We conclude this section by considering all the intervals and all
the curves of indices $0$ and $1$ with respect to the fixed
fundamental system $\alpha_1, \dots, \alpha_n$.

The intervals of index $0$ are the $\widehat x_{i,j}$'s. The
curves of index $0$ are the curves $\alpha_{i,j}$ with $1 \leq
i,j \leq n$, defined in the following way: $\alpha_{i,i} =
\alpha_i$, $\alpha_{i,j} = (\alpha_i)\widehat x_{i,j}^{-1}$ if $i
< j$, $\alpha_{i,j} = (\alpha_i)\widehat x_{i,j}$ if $j < i$. Such
intervals and curves are related by the following equalities:
$$(\alpha_{i,j})\widehat x_{j,k}^{-1} = \cases{
   \hbox to 28pt{$\alpha_{i,k}$}
      \ \ \hbox{if}\ \ i \leq j < k\ \ \hbox{or}\ \ j < k \leq i\,,\cr
   \vrule width 0pt height 12pt
     \hbox to 28pt{$\alpha_{i-1,k}$}\ \ \hbox{if}\ \ k < i \leq j\,;}$$
$$(\alpha_{i,j})\widehat x_{j,k} = \cases{
   \hbox to 28pt{$\alpha_{i,k}$\hfill}
      \ \ \hbox{if}\ \ i \leq k < j\ \ \hbox{or}\ \ k < j \leq i\,, \cr
   \vrule width 0pt height 12pt
      \hbox to 28pt{$\alpha_{i+1,k}$\hfill}\ \ \hbox{if}\ \ j \leq i < k\,.}$$

\medskip

The intervals of index $1$ are the intervals $\widehat x_{i,j,k}$
with $1 \leq i,j,k \leq n$ such that $i < k$ and $i \neq j \neq
k$, given by: $\widehat x_{i,j,k} = (\widehat x_{i,j})\widehat
x_{j,k}$ if $i < j < k$ and $\widehat x_{i,j,k} = (\widehat
x_{i,j})\widehat x_{j,k}^{-1}$ if $i < k < j$ or $j < i < k$. As
a notational convenience, we also set $\widehat x_{i,j,k} =
\widehat x_{k,j,i}$ if $i > k$ and $\widehat x_{i,i,j} = \widehat
x_{i,j,j} = \widehat x_{i,j}$ for every $i \neq j$.

Finally, the curves of index $1$ are the curves $\alpha_{i,j,k}$
with $1 \leq i,j,k \leq n$ such that $i \neq j \neq k$, defined
as follows:
$$\alpha_{i,j,k} = \cases{
(\alpha_{i,j})\widehat x_{j,k}
   \ \ \hbox{if}\ \ i<j<k\ \ \hbox{or}\ \ j<k \leq i\ \ \hbox{or}\ \ k<i<j\,,
\cr\vrule width 0pt height14pt (\alpha_{i,j})\widehat x_{j,k}^{-1}
   \ \ \hbox{if}\ \ k<j<i\ \ \hbox{or}\ \ j<i<k\ \ \hbox{or}\ \ i \leq k<j\,.}$$

\section{Liftable braids with respect to $p:B^2 \to
B^2$\label{specialcase/sec}}

\indent\indent By the results of Section \ref{coverings/sec}, for
every $n \geq 1$ there exists a unique (up to equivalence) simple
branched covering $p_n: B^2 \to B^2$ of order $d = n+1$ with $n$
branch points. Moreover, the $p_n$'s represent (up to equivalence)
all the coverings of $B^2$ onto itself.

We assume the base point $* \in S^1$, the branch points $P_1,
\dots, P_n \in \Int B^2$, the fundamental system $\alpha_1,
\dots, \alpha_n$ and the numbering of the sheets of $p_n$ fixed
in such a way that: (1) $\alpha_i$ joins $*$ to $P_i$ for every
$i=1, \dots, n$; (2) the monodromy sequence $\phi(\alpha_1),
\dots, \phi(\alpha_n)$ is in the canonical form $(1\;2), \dots,
(d{-}1\;d)$ given in the proof of Theorem A, namely
$\phi(\alpha_i) = (i\;i{+}1)$ for every $i = 1, \dots, n$.

\medskip

{\sl In this section, all the curves $\alpha_{i,j}$ and
$\alpha_{i,j,k}$ and all the intervals $x_i$, $x_{i,j}$,
$x_{i,j,k}$, $\widehat x_{i,j}$, $\widehat x_{i,j,k}$, as well as
all the indexes of curves and intervals, except where expressly
indicated, are referred to the fundamental system  $\alpha_1,
\dots, \alpha_n$.}

\medskip

In order to prove Theorem C, let us begin with some preliminary
results about curves. We recall that $\L_n \subset \B_n$ denotes
the subgroup of the liftable braids with respect to $p_n$.

\medskip

By direct computation we get the following monodromies:
$$\phi(\alpha_{i,j}) =
   \cases{(i\;j{+}1)\ \ \ \hbox{if}\ \ i \leq j\,,\cr
          (i{+}1\;j)\ \ \ \hbox{if}\ \ j \leq i\,;}$$
$$\phi(\alpha_{i,j,k}) = \cases{
   \hbox to 50pt{$(j{+}1\;k{+}1)$\hfill}\ \ \ \hbox{if}\ \ i < j < k \
      \ \hbox{or}\ \ i \leq k < j\,,\cr
   \hbox to 50pt{$(j\;k)$\hfill} \ \ \
   \hbox{if}\ \ k < j < i \ \ \hbox{or}\ \ j < k \leq i\,,\cr
   \hbox to 50pt{$(j{+}1\;k)$\hfill}\ \ \ \hbox{if}\ \ k \leq i < j\,,\cr
   \hbox to 50pt{$(j\;k{+}1)$\hfill}\ \ \ \hbox{if}\ \ j < i \leq k\,.}$$

\medskip

Assuming $n > 1$, we say that a curve $\alpha$ for $p_n$ is {\sl
regular} if $p_n^\alpha$ is equivalent to $p_{n-1} \sqcup
\id_{B^2}$. We observe that, if $\alpha$ is a regular curve then
$(\alpha)b$ is a regular curve for any liftable braid $b \in
\L_n$.

\medskip \begin{lemma} \label{regularcurves/thm}% Lemma 2.1.
The curve $\alpha_j$ is regular for every $j = 1,\dots, n$.
Moreover, the equivalence between $p_{n-1}$ and the non-trivial
component of $p_n^{\alpha_j}$ is induced by a homeomorphism
$h_j:B^2 \to B^2_{\alpha_j}$ such that: $h_j(\alpha_i) =
\alpha_i'$ for $1 \leq i < j$; $h_j(\alpha_i) = \alpha_{i+1}'$
for $j \leq i \leq n-1$; $h_j(x_i) = x_i$ for $1 \leq i < j-1$;
$h_j(x_{j-1}) = (x_{j-1})x_j^{-1}$; $h_j(x_i) = x_{i+1}$ for $j-1
< i \leq n-2$.
\end{lemma}

\medskip\begin{proof}
The fundamental system $\alpha'_1, \dots, \alpha'_{j-1},
\alpha'_{j+1}, \dots, \alpha'_n$ for $p_n^{\alpha_j}$ has
monodromy sequence $(1\;2), \dots, (j{-}2\;j{-}1),
(j{-}1\;j{+}1), (j{+}1\;j{+}2), \dots, (n{-}1\;n)$. Let $h_j:B^2
\to B^2_{\alpha_j}$ be the homeomorphism uniquely determined (up
to isotopy) by $h_j(\alpha_i) = \alpha'_i$ for $1 \leq i < j$ and
$h_j(\alpha_i) = \alpha'_{i+1}$ for $j \leq i \leq n-1$. By the
Lifting theorem, $h_j$ can be lifted to give an equivalence
between $p_{n-1} \sqcup{\id_{B^2}}$ and $p_n^{\alpha_j}$. Hence,
$h_j$ induces an equivalence between $p_{n-1}$ and the
non-trivial component of $p_n^{\alpha_j}$. A straightforward
computation of the intervals $h_j(x_i)$ completes the proof.
\end{proof}

\medskip \begin{lemma} \label{regularcurves1n/thm} % Lemma 2.2.
The curves $\alpha_{1,n}$ e $\alpha_{n,1}$ are regular. Moreover,
we have that: the equivalence between $p_{n-1}$ and the
non-trivial component of $p_n^{\alpha_{1,n}}$ is induced by a
homeomorphism $h_{1,n}:B^2 \to B^2_{\alpha_{1,n}}$ such that
$h_{1,n}(\alpha_i) = \alpha_i''$ for $1 \leq i \leq n-1$ and
$h_{1,n}(x_i) = x_i$ for $1 \leq i \leq n-2$; the equivalence
between $p_{n-1}$ and the non-trivial component of
$p_n^{\alpha_{n,1}}$ is induced by a homeomorphism $h_{n,1}:B^2
\to B^2_{\alpha_{n,1}}$ such that $h_{n,1}(\alpha_i) =
\alpha_{i+1}'$ for $1 \leq i \leq n-1$ and $h_{n,1}(x_i) =
x_{i+1}$ for $1 \leq i \leq n-2$.
\end{lemma}

\medskip\begin{proof}
Similar to the previous one, except that we consider the
fundamental system $\alpha_1'', \dots, \alpha_{n-1}''$ instead of
$\alpha_1', \dots, \alpha_{n-1}'$ for the covering
$p_n^{\alpha_{1,n}}$.
\end{proof}

\medskip \begin{lemma} \label{regularcurves0/thm} % Lemma 2.3.
The only  regular curves of index $0$ are $\alpha_1, \dots,
\alpha_n$, $\alpha_{1,n}$ and $\alpha_{n,1}$. Among these, only
$\alpha_{1,n}$ and $\alpha_{n,1}$ are $\L_n$-equivalent to each
other.
\end{lemma}

\medskip\begin{proof}
Lemmas \ref{regularcurves/thm} and \ref{regularcurves1n/thm} say
that the curves $\alpha_1, \dots,\alpha_n, \alpha_{1,n}$ and
$\alpha_{n,1}$ are regular. In the previous section, we observed
that any other curve of index $0$ have to be an $\alpha_{i,j}$
with $j \neq i$ and $(1,n) \neq (i,j) \neq (n,1)$. If $j > i$,
then the curves $\alpha_1, \dots, \alpha_{i-1}, \alpha_{i,j},
\alpha_{i}, \dots, \alpha_{j-1}, \alpha_{j+1}, \dots, \alpha_n$
constitute a fundamental system for $p_n$ with monodromy sequence
   $$(1\;2), \dots, (i{-}1\;i), (i\;j{+}1), (i\;i{+}1), \dots, (j{-}1\;j),
     (j{+}1\;j{+}2), \dots, (n\;n{+}1).$$
  If $j < i$ then the curves $\alpha_1, \dots ,\alpha_{j-1},\alpha_{j+1}, \dots,
\alpha_i, \alpha_{i,j},\alpha_{i+1}, \dots, \alpha_n$ constitute
a fundamental system for $p_n$ with sequence of monodromies
   $$(1\;2), \dots, (j{-}1\;j), (j{+}1\;j{+}2), \dots, (i\;i{+}1), (j\;i{+}1),
     (i{+}1\;i{+}2), \dots,(n\;n{+}1).$$
In both cases, none of the curves $\alpha_{i,j}$ is regular, as
can be immediately proved by using Lemma \ref{disk/thm}.

For the second part of the lemma, we observe that the monodromies
of the curves taken into account are distinct from each other,
with the only exception of $\phi(\alpha_{1,n}) =
\phi(\alpha_{n,1}) = (1\;n{+}1)$. On the other hand, since
$\alpha_{n,1} = (\alpha_{1,n})b$, with $b = (x_{n-1}\dots
x_1)^{n+1} \in \L_n$, we have that $\alpha_{1,n}$ and
$\alpha_{n,1}$ are $\L_n$-equivalent.
\end{proof}

\medskip \begin{lemma} \label{tworegularcurves/thm} % Lemma 2.4.
Any fundamental system $\beta_1, \dots, \beta_n$ for $p_n$ with
$n > 1$, contains at least two regular curves $\beta_{i_1}$ e
$\beta_{i_2}$.
\end{lemma}

\medskip\begin{proof}
Let $\Gamma = \Gamma_{p_n}(\beta_1, \dots, \beta_n)$ be the graph
associated to $\beta_1, \dots, \beta_n$. Moreover, for any $i =
1, \dots, n$, let $\Gamma_i = \Gamma_{p_n^{\beta_i}}(\beta_1',
\dots, \beta_{i-1}', \beta_{i+1}', \dots, \beta_n')$ be the graph
associated to the fundamental system $\beta_1', \dots,
\beta_{i-1}', \beta_{i+1}', \dots, \beta_n'$ for $p_n^{\beta_i}$.
By Lemma \ref{homotopy/thm}, $\Gamma$ is a tree. On the other
hand, it follows from Lemma \ref{disk/thm} that all the
$\Gamma_i$'s have two connected components and that $\beta_i$ is
regular if and only if one component of $\Gamma_i$ consists of a
single vertex. Then, it is enough to prove that there exist two
graph $\Gamma_{i_1}$ and $\Gamma_{i_2}$ with that property.

The graph $\Gamma_i$ can obtained from $\Gamma$, by removing the
edge $e_i$ and replacing the edge $e_l = \{v_{j_l},v_{k_l}\}$
with the new edge $e_{l-1} =
\{v_{\phi(\beta_i)(j_l)},v_{\phi(\beta_i)(k_l)}\}$, for every $l >
i$. We remark that the edges $e_1, \dots, e_{i-1}$, as well as
all the $e_l$'s not meeting $e_i$, are left unaltered.

Now let $\Gamma'$ be the full subgraph of $\Gamma$ generated by
all the vertices of valence greater than $1$. It is not difficult
to see that $\Gamma$ collapses to $\Gamma'$ (remember that $n >
1$). Then, also $\Gamma'$ is a non-empty tree.

If $\Gamma'$ reduces to a single vertex, this vertex is contained
in all the edges $e_1, \dots, e_n$ of $\Gamma$. In this case, we
have that $\Gamma_1$ and $\Gamma_n$ have the required property.
Otherwise, $\Gamma'$ must contain al least two different valence
one vertices $w_1$ and $w_2$. From these vertices come out two
different edges $e_{i_1}$ and $e_{i_2}$ of $\Gamma - \Gamma'$,
such that the graphs $\Gamma_{i_1}$ e $\Gamma_{i_2}$ have the
required property.

Let us see how to determine $i_1$ (in the same way could be
determined $i_2$). Let $e_{l_1}$ be the only edge of $\Gamma'$
containing $w_1$. Since the valence of $w_1$ in $\Gamma$ is
greater than one, there is least one edge of $\Gamma - \Gamma'$
containing $w_1$. Then, we can set $i_1$ equal to the maximum
among the indices of such edges.
\end{proof}

We continue by considering some properties of the intervals.
First of all, we observe that all the intervals $x_i$ are of type
$3$ with respect to  $p_n$, while all the intervals $x_{i,j}$
with $j > i +1$ are of type $2$.

\medskip \begin{lemma} \label{type3/thm} % Lemma 2.5.
All the index\/ $0$ intervals are of type $3$ with respect to
$p_n$.
\end{lemma}

\medskip\begin{proof}
We recall that the index $0$ intervals are the $\widehat
x_{i,j}$'s with $i < j$. Such intervals are of type $3$ by Lemma
\ref{intervals/thm}, since the curve $\alpha_i$ meets $\widehat
x_{i,j}$ only at its endpoint, $\phi(\alpha_i) = (i\;i{+}1)$ and
$\phi((\alpha_i)\widehat x_{i,j}) = \phi(\alpha_{i+1,j}) =
(i{+}1\;j{+}1)$.
\end{proof}

\medskip \begin{lemma} \label{type2/thm} % Lemma 2.6.
All the index\/ $1$ intervals are of type $2$ with respect to
$p_n$.
\end{lemma}

\medskip\begin{proof}
We recall that the index $1$ intervals are the $\widehat
x_{i,j,k}$'s with $i < k$ and $i \neq j \neq k$. The curve
$\alpha_i$ meets $\widehat x_{i,j,k}$ only at its endpoint and we
have that $(\alpha_i)\widehat x_{i,j,k}$ coincides with
$\alpha_{i+1,j,k}$ if $i < j < k$ or $i < j < k$ and with
$\alpha_{i,j,k}$ if $j < i < k$. In any case, the transpositions
$\phi(\alpha_i)$ and $\phi((\alpha_i)\widehat x_{i,j,k})$ are
disjoint. Then, $\widehat x_{i,j,k}$ is of type $2$ by Lemma
\ref{intervals/thm}.
\end{proof}

\medskip \begin{lemma} \label{type1/thm} % Lemma 2.7.
There are no intervals of type $1$ with respect to $p_n$.
\end{lemma}

\medskip\begin{proof}
Given any interval $x$ and any curve $\alpha$ which meets $x$
only at its endpoint, let $\beta_1, \dots, \beta_n$ be any
fundamental system such that $\beta_1 = \alpha$ and $\beta_2 =
(\alpha)x$. If $x$ were of type $1$, $\phi(\beta_1)$ would
coincide with $\phi(\beta_2)$, in contradiction with Lemma
\ref{homotopy/thm} and Lemma \ref{disk/thm}.
\end{proof}

For sake of simplicity, we denote by $\I_n \subset \L_n$ the
group  $\I_{p_n}$ generated by the liftable powers of intervals.
The braids $x_i^3$ and $x_{i,j}^2$ with $1 \leq i < n$ and $i+1 <
j \leq n$ belong to $\I_n$. In fact, we will see that they
generate $\I_n$.

\medskip \begin{lemma} \label{specialcurves0/thm}% Lemma 2.8.
If $\alpha$ is a curve whose interior meets each one of the curves
$\alpha_1, \dots, \alpha_n$ in at most one point, then $\alpha$
is $\I_n$-equivalent to a curve of index $0$.
\end{lemma}

\medskip\begin{proof}
We proceed by induction on the index of $\alpha$, assuming that
$\alpha$ minimizes the number of intersection points with
$\alpha_1 \cup \dots \cup \alpha_n$ in its isotopy class.

We start with the index 1 case. In this case, we have the curves
$\alpha = \alpha_{i,j,k}$, with $1 \leq i,j,k \leq n$ such that
$i \neq j \neq k$, defined in Section \ref{braids/sec}. If $i =
k$, it suffices to observe that $\alpha_{i,j,i}$ is
$\I_n$-equivalent to the index $0$ curve $(\alpha_{i,j,i})
\widehat x_{i,j}^{\pm3} = \alpha_{i\pm1,j}$, where $\pm$ is the
sign of $j - i$, being $\widehat x_{i,j}$ of type $3$ by Lemma
\ref{type3/thm}. If $i \neq k$, then $\alpha_{i,j,k}$ is
$\I_n$-equivalent to $(\alpha_{i,j,k}) \widehat x_{i,j,k}^{\pm2}$,
where $\pm$ is the sign of $j - i$, being $\widehat x_{i,j,k}$ of
type $2$ by Lemma \ref{type2/thm}. The curve $(\alpha_{i,j,k})
\widehat x_{i,j,k}^{\pm2}$ has index $0$ if $|i - j| = 1$, while
it coincides with the curve $\alpha_{i\pm1,j,k}$, if $|i - j| >
1$. So, we can conclude the case of the $\alpha_{i,j,k}$'s with
$i \neq k$, by induction on $|i - j| \geq 1$.

Now we suppose that $\alpha$ has index ${} > 1$. Let $P_k$ be the
endpoint of $\alpha$ and let $Q_i \in \alpha \cap \alpha_i$ and
$Q_j \in \alpha \cap \alpha_j$ be respectively the last but one
and the last point in which the interior of $\alpha$ (oriented
from $*$ to $P_k$) meets the  curves $\alpha_1, \dots, \alpha_n$.
We consider the following  arcs: $t_i \subset \alpha_i$ with
endpoints $Q_i$ and $P_i$, $t_j \subset \alpha_j$ with endpoints
$Q_j$ and $P_j$, $s_i \subset \alpha$ with endpoints $Q_i$ and
$P_k$, $s_j \subset \alpha$ with endpoints $Q_j$ and $P_k$. By
hypothesis we have $i \neq j$. Moreover, we can assume $j \neq
k$, otherwise we could remove the intersection $Q_j$ up to
isotopy.

If $i = k$, the interval $x = t_j \cup s_j$ has index $0$. Then,
by Lemma \ref{type3/thm}, $\alpha$ is $\I_n$-equivalent to the
curve $(\alpha)x^{\pm3}$, with sign $-$ if $t_j$ is on the left
of $\alpha$ and sign $+$ if $t_j$ is on the right of $\alpha$.
The curve $(\alpha)x^{\pm3}$ has index less than $\alpha$ (the
intersections $Q_i$ and $Q_j$ disappear) and it is
$\I_n$-equivalent to a  curve of index  $0$ by the induction
hypothesis.

If $i \neq k$, the interval $x = t_i \cup s_i$ has index $1$.
Then, by Lemma \ref{type2/thm}, $\alpha$ is $\I_n$-equivalent to
the curve $(\alpha)x^{\pm2}$, with sign $-$ if $t_i$ is on the
left of $\alpha$ and sign $+$ if $t_i$ is on the right of
$\alpha$. The curve $(\alpha)x^{\pm2}$ has index less than
$\alpha$ (the intersection $Q_i$ disappears) and it is
$\I_n$-equivalent to a curve of index $0$ by the induction
hypothesis.
\end{proof}

\medskip \begin{lemma} \label{curves0/thm} % Lemma 2.9.
Every curve $\alpha$ is $\I_n$-equivalent to a curve of index $0$.
\end{lemma}

\medskip\begin{proof}
We proceed by induction on $n$. For $n = 1$ there is nothing to
prove. So, let us suppose $n > 1$. First of all, we consider the
special case in which $\alpha \cap \alpha_j = \{*\}$ for some $j
= 1, \dots, n$.  By Lemma \ref{regularcurves/thm} and by the
induction hypothesis, it exists a braid $b \in
\I_{p_n^{\alpha_{\smash{j}}}}$ such that the curve $(\alpha')b$
has index $0$ with respect to the fundamental system $\alpha'_1,
\dots, \alpha'_{j-1}, \alpha'_{j+1},\dots,\alpha'_n$ for
$p_n^{\alpha_j}$. The braid $b$ can also be considered as a braid
in $\I_n$ and it is easy to verify that the curve $(\alpha)b$
satisfies Lemma \ref{specialcurves0/thm}. Then $\alpha$ is
$\I_n$-equivalent to a curve of index $0$. By Lemma
\ref{regularcurves1n/thm}, also the cases $\alpha \cap
\alpha_{1,n} = \{*\}$ and $\alpha \cap \alpha_{n,1} = \{*\}$, with
the braid $b$ respectively in $\I_{p_n^{\alpha_{\smash{1,n}}}}$
and in $\I_{p_n^{\alpha_{\smash{n,1}}}}$ can be treated in an
analogous way.

Now we carry on the proof by induction on the index of $\alpha$,
assuming that $\alpha$ meets every $\alpha_j$ in some point other
than $*$. For every $j = 1, \dots, n$, we denote by $Q_j$ the
point of $\alpha \cap \alpha_j$ nearest to $P_j$ along
$\alpha_j$, and by $\beta_j$ the curve obtained following $\alpha$
from $*$ to $Q_j$ and then $\alpha_j$ from $Q_j$ to $P_j$. If
$P_k$ is the endpoint of $\alpha$, then $\beta_k = \alpha$ and
all the curves $\beta_j$ with $j \neq k$ have index less than
$\alpha$. Since the curves $\beta_1,\dots,\beta_n$, suitably
renumbered, constitute a fundamental system, Lemma
\ref{tworegularcurves/thm} ensures the existence of $l \neq k$
such that $\beta_l$ is regular. By the induction hypothesis, there
exists $b \in \I_n$ such that $(\beta_l)b$ has index $0$. Then
$(\beta_l)b$ coincides either with some $\alpha_j$ or with
$\alpha_{1,n}$ or with $\alpha_{n,1}$. Hence, $(\alpha)b$ is
$\I_n$-equivalent to a curve of index $0$, being included in the
cases examined at the beginning of the proof. It follows that
$\alpha$ is as well $\I_n$-equivalent to a curve of index $0$.
\end{proof}

\medskip \begin{lemma} \label{liftablepowers/thm} % Theorem 2.10.
$\L_n$ is generated by liftable powers of intervals.
\end{lemma}

\medskip\begin{proof}
We proceed by induction on $n$. If $n = 1$ there is nothing to
prove. If $n > 1$ and $b \in \L_n$, then Lemmas \ref{curves0/thm}
and \ref{regularcurves0/thm}, give us a braid $c \in \I_n$ such
that $(\alpha_n)bc = \alpha_n$, in such a way that $bc$ can be
considered as a braid in $\L_{p_n^{\alpha_{\smash{n}}}}$. By the
regularity of $\alpha_n$ and by induction hypothesis, we have $bc
\in \I_{p_n^{\alpha_{\smash{n}}}} \subset \I_n$ and therefore $b
\in \I_n$.
\end{proof}

Now,  let $\J_n \subset \L_n$ denote the subgroup generated by the
braids $x_i^3$ and $x_{i,j}^2$ with $1 \leq i < n$ and $i+1 < j
\leq n$. We want to prove that actually $\J_n = \L_n$, that is
our Theorem C.

To get this goal, observe that in the proof of Lemma
\ref{liftablepowers/thm} we do not use  the liftable powers of
all the intervals, but only of some particular intervals.
Therefore, it is enough to show that each one of these particular
intervals is $\J_n$-equivalent to some $x_i$ or $x_{i,j}$.

\medskip \begin{lemma} \label{specialliftablepowers} % Lemma 2.11.
Every interval $x=(x_i)x_{i+1}^{e_{i+1}}\dots x_{j-1}^{e_{j-1}}$,
with $e_{i+1}, \dots, e_{j-1} = \pm 1$ and $1 \leq i < j \leq n$,
is $\J_n$-equivalent to some $x_{h,k}$, so all the liftable
powers of $x$ belong to $\J_n$.
\end{lemma}

\medskip\begin{proof}
By induction on the number of negative $e_l$'s. If all the $e_l$'s
are positive, then $x = x_{i,j}$. Otherwise, let $m \geq i+1$ be
the minimum integer such that $e_m = -1$. If $m=i+1$, then $x =
(x_i)x_{i+1}^{-1}x_{i+2}^{e_{i+2}}\dots x_{j-1}^{e_{j-1}} =
(y)z^2x_i^3$ with $y = (x_{i+1})x_{i+2}^{e_{i+2}}\dots
x_{j-1}^{e_{j-1}}$ and $z = (x_i)x_{i+1}x_{i+2}^{e_{i+2}}\dots
x_{j-1}^{e_{j-1}}$. Since $y$ and $z$ are $\J_n$-equivalent to
some $x_{h,k}$ by the induction hypothesis and $z$ is of type $2$
(so $z^2 \in \J_n$), we have that also $x$ is $\J_n$-equivalent to
some $x_{h,k}$. If $m > i+1$, then $x = (x_i)x_{i+1} \dots
x_{m-1}x_m^{-1}x_{m+1}^{e_{m+1}}\dots x_{j-1}^{e_{j-1}} =
(t)x_{i,m}^2$ with $t = (x_i)x_{i+1} \dots
x_{m-1}x_mx_{m+1}^{e_{m+1}}\\ \dots x_{j-1}^{e_{j-1}}$. Since $t$
is $\J_n$-equivalent to some $x_{h,k}$ by the induction
hypothesis, also $x$ is $\J_n$-equivalent to some $x_{h,k}$.
\end{proof}

\medskip \begin{lemma} \label{liftablepowers1} % Lemma 2.12.
Every interval $x$ of index ${} \leq 1$ is $\J_n$-equivalent to
some $x_{h,k}$, so all the liftable powers of $x$ belong to
$\J_n$.
\end{lemma}

\medskip\begin{proof} The intervals of index $0$, that is the $\widehat
x_{i,j}$'s have been already considered in Lemma
\ref{specialliftablepowers}. The same also holds for the intervals
of index $1$ of type $\widehat x_{i,j,k}$ with $i < j < k$, in
fact for these intervals we have $\widehat x_{i,j,k} =
(x_i)x_{i+1}^{-1} \dots x_{j-1}^{-1} x_j x_{j+1}^{-1} \dots \\
x_{k-1}^{-1}$.

It remains only to deal with the intervals $\widehat x_{i,j,k} =
(\widehat x_{i,j}) \widehat x_{j,k}^{-1}$ such that either $i < k
< j$ or $j < i < k$. In the first case we have that $\widehat
x_{i,j,k}$ is $\J_n$-equivalent to the interval $(\widehat
x_{i,j,k}) \widehat x_{j,k}^3 = \widehat x_{i,k,j}$. In the second
case we have that $\widehat x_{i,j,k}$ is $\J_n$-equivalent to
the interval $(\widehat x_{i,j,k}) \widehat x_{i,j}^{-3} =
\widehat x_{j,i,k}$. Hence, in both the cases $\widehat
x_{i,j,k}$ is $\J_n$-equivalent to an interval having the form
considered above.
\end{proof}

\medskip\begin{proof}[Proof of Theorem C]
We proceed by induction on $n$. For $n = 1$ there is nothing to
prove. So, let us suppose $n > 1$. In the proof of Lemma
\ref{specialcurves0/thm}, the $\I_n$-equivalence desired is
obtained by using liftable powers of intervals of index ${} \leq
1$, which belong in $\J_n$ by Lemma \ref{liftablepowers1}. On the
other hand, in proofs of Lemmas \ref{curves0/thm} and
\ref{liftablepowers/thm}, we use liftable powers of intervals in
$\I_{p_n^{\alpha_{\smash{j}}}}$, $\I_{p_n^{\alpha_{\smash{1,n}}}}$
and $\I_{p_n^{\alpha_{\smash{n,1}}}}$. By the induction
hypothesis, these groups are generated by braids of the form
$y_i^3$ and $y_{h,k}^2$ with $y_i = h(x_i)$ and $y_{h,k} =
h(x_{h,k})$, where $h$ denotes one of the homeomorphism $h_j$,
$h_{1,n}$ and $h_{n,1}$ given by Lemmas \ref{regularcurves/thm}
and \ref{regularcurves1n/thm}. It is not difficult to see that the
intervals $y_i$ and $y_{h,k}$ are among the ones considered in
Lemma \ref{specialliftablepowers}, so their liftable powers belong
to $\J_n$.

Then, we can replace the group $\I_n$ with the group $\J_n$ in
Lemmas \ref{specialcurves0/thm} and \ref{curves0/thm} as well as
in the proof of Lemma \ref{liftablepowers/thm}, in order to get
$\L_n = \J_n$.
\end{proof}

\section{Liftable braids with respect to $p:F \to B^2$\label{generalcase/sec}}

\indent\indent All this section is devoted to prove Theorem B.
Here, we consider an arbitrary connected simple branched covering
$p:F \to B^2$ of order $d$ with $n$ branch points. As in the
previous section, we assume the base point $* \in S^1$, the branch
points $P_1, \dots, P_n \in \Int B^2$, the fundamental system
$\alpha_1, \dots, \alpha_n$ and the numbering of the sheets of $p$
fixed in such a way that: (1) $\alpha_i$ joins $*$ to $P_i$ for
every $i = 1, \dots, n$; (2) the monodromy sequence
$\phi(\alpha_1), \dots, \phi(\alpha_n)$ is in the canonical form
given in the proof of Theorem A.

\medskip \begin{lemma} \label{betagamma/thm} % Lemma 3.1.
Let $\beta$ be a curve such that $p^\beta$ is connected and let
$\beta_1, \dots, \beta_n$ be a fundamental system for $p$. Then
$\beta$ is $\I_p$-equivalent to a curve $\gamma$ such that
$\gamma \cap \beta_i = \{*\}$ for some $i = 1, \dots, n$.
\end{lemma}

\medskip\begin{proof}
Let $\gamma$ be a curve of minimum index with respect to the
fundamental system $\beta_1, \dots, \beta_n$ among all the curves
$\I_p$-equivalent to $\beta$. Let us also assume that $\gamma$
minimizes the number of intersection points with $\beta_1 \cup
\dots \cup \beta_n$ in its isotopy class. We claim that there
exists an integer $i = 1, \dots, n$ such that $\gamma \cap
\beta_i = \{*\}$.

Suppose, by the contrary, that $\gamma$ meets any $\beta_i$ in
some point other than $*$. For each $i=1,\dots,n$, we denote by
$Q_i$ the last point of $\gamma\cap\beta_i$ along $\beta_i$
(starting from $*$) and with $\gamma_i$ the curve obtained
following $\gamma$ until $Q_i$ and then $\beta_i$ until its
endpoint. Up to isotopy, we can suppose
$\gamma_i\cap\gamma_j=\{*\}$ for all $i \neq j$. If the endpoint
of $\gamma$ coincides with the endpoint of $\beta_k$, then
$\gamma_k = \gamma$ and any curve $\gamma_i$ with $i \neq k$ has
index less than $\gamma$. We denote by $\sigma_i =
\phi(\gamma_i)$ the monodromy of $\gamma_i$. In particular, let
$\sigma_k = (a\;b)$ be the monodromy  of $\gamma$.

Let us consider the intervals $y_{i,j} \simeq \gamma_i \cup
\gamma_j$ for $i \neq j$ and $1 \leq i,j \leq n$. We observe that
all the $y_{i,k}$'s are of type $3$, that is any transposition
$\sigma_i$ with $i \neq k$ is distinct but not disjoint from
$(a\;b)$. Indeed, if $y_{i,k}$ were of type $1$ or $2$ then
$\gamma$ would be $\I_p$-equivalent to the curve
$(\gamma)y_{i,k}^{\pm2}$, with $-$ or $+$ depending on whether
$\gamma_i$ is on the left or on the right of $\gamma$, which has
index less than $\gamma$.

On the other hand, if $\gamma_i$ and $\gamma_j$, with $i,j \neq
k$, are on the same side with respect to $\gamma$, then
$\{\sigma_i, \sigma_j\} \neq \{(a\;b),(b\;c)\}$. Indeed, assuming
that $Q_i$ precedes $Q_j$ along $\gamma$ (starting from $*$), the
equality $\{\sigma_i, \sigma_j\} = \{(a\;b),(b\;c)\}$ would imply
the liftability of the interval $x = (y_{j,k})y_{i,j}^{\pm2}$,
with $-$ or $+$ depending on the fact that $\gamma_i$ and
$\gamma_j$ are on the left or on the right of $\gamma$. Therefore,
$\gamma$ would be $\I_p$-equivalent to the curve $\delta =
(\gamma)x^{\pm1}$, with the same choice for the sign, which has
index less than $\gamma$.

Analogously, if $\gamma_i$ and $\gamma_j$, with $i,j \neq k$, are
on opposite sides with respect to $\gamma$, then $\sigma_i \neq
\sigma_j$. Indeed, assuming as above that $Q_i$ precedes $Q_j$
along $\gamma$ (starting from $*$), the equality $\sigma_i =
\sigma_j$ would imply the liftability of $y_{i,j}$. Therefore,
$\gamma$ would be $\I_p$-equivalent to the curve $\delta =
(\gamma)y_{i,j}^{\pm1}$, with $-$ or $+$ depending on the fact
that $\gamma_i$ is on the left or on the right of $\gamma$, which
has index less than $\gamma$.

Hence, by renumbering the $\gamma_i$'s in clockwise order, we get
a new fundamental system for $p$, whose monodromy sequence has
the form $$(c_1\;d_1), \dots, (c_{h-1}\;d_{h-1}), (a\;b),
(c_{h+1}\;d_{h+1}), \dots, (c_n\;d_n)$$ and satisfies the
following properties: $c_i\not\in\{a,b\}$ and $d_i\in \{a,b\}$
for any $i \neq h$;  if $i,j < h$ or $i,j > h$ then $c_i = c_j
\Rightarrow d_i = d_j$; if $i < h < j$ then $c_i = c_j
\Rightarrow d_i \neq d_j$. Then, by putting $C_a^- =
\{c_i\;|\;d_i = a \wedge i < h\}$, $C_a^+ = \{c_i\;|\;d_i = a
\wedge i > h\}$, $C_b^- = \{c_i\;|\;d_i = b \wedge i < h\}$ and
$C_b^+ = \{c_i\;|\;d_i = b \wedge i > h\}$, we have $C_a^- \cap
C_b^- = C_a^+ \cap C_b^+ = C_a^- \cap C_a^+ = C_b^- \cap C_b^+ =
\emptyset$.

Now, the fundamental system $\gamma_1', \dots,
\gamma_{h-1}',\gamma_{h+1}', \dots, \gamma_n'$ for the covering
$p^\gamma$ has monodromy sequence $(c_1\;d_1), \dots,
(c_{h-1}\;d_{h-1}), (c_{h+1}\;\bar d_{h+1}), \dots, (c_n\;\bar
d_n)$, where $\bar d_i = a$ if $d_i= b$ and $\bar d_i = b$ if $d_i
= a$. Such a sequence of transpositions can be reordered in the
form $(e_1\;a), \dots, (e_l\;a), (e_{l+1}\;b), \dots,
(e_{n-1}\;b)$ with $e_i \in C_a^- \cup C_b^+$ if $i \leq l$ and
$e_i \in C_a^+ \cup C_b^-$ if $i \geq l+1$. Therefore the two sets
$C_a^- \cup C_b^+ \cup \{a\}$ e $C_a^+ \cup C_b^- \cup \{b\}$ are
disjoint, non-empty and closed with respect to the action of the
group generated by these transpositions. Of course, this fact
contradicts the connection of $p^\gamma \cong p^\beta$. So,
$\gamma$ cannot meet any $\beta_i$ in some point other than the
point $*$.
\end{proof}

\medskip \begin{lemma} \label{betadelta/thm} % Lemma 3.2.
Let $\beta$ be a curve such that $\beta = (\alpha_m)b$, with $b
\in \L_p$ and $1 \leq m \leq n$, and $p^\beta \cong p^{\alpha_m}$
is connected. Then $\beta$ is $\I_p$-equivalent to a curve
$\delta$ such that $\delta \cap \alpha_i = \{*\}$ for some $i =
1, \dots, m$ and $\delta$ starts from $*$ on the left (resp.
right) of $\alpha_i$ if  $i < m$ (resp. $i \geq m$).
\end{lemma}

\medskip\begin{proof}
By Lemma \ref{betagamma/thm}, $\beta$ is $\I_p$-equivalent to a
curve $\gamma$ which meets at least one of the $\alpha_i$'s only
in $*$. In other words, the set $S \subset \{1,\dots,n\}$ of the
$i$'s such that $\gamma \cap \alpha_i = \{*\}$ is nonempty. We
can also assume that $\gamma$ has minimum index (with respect to
the fundamental system $\alpha_1, \dots, \alpha_n$) among all the
curves having such property in the $\I_p$-equivalence class of
$\beta$.

If there exists $i \in S$ such that either $i < m$ and $\gamma$
starts from $*$ on the left of $\alpha_i$ or $i \geq m$ and
$\gamma$ starts from $*$ on the right of $\alpha_i$, then we can
put $\delta = \gamma$.

If such an $i$ does not exist, but there exists $i \in S$ such
that the interval $x \simeq \gamma \cup \alpha_i$ is of type $1$
or $2$, then we can put $\delta = (\gamma) x^{\pm2}$, with $+$ or
$-$ depending on the fact that $\gamma$ starts from $*$ on the
left or on the right of $\alpha_i$.

%%%%%%%%%%%%

In the remaining cases, all the curves $\alpha_i$ with $i \in S$
have the same monodromy and start from $*$ on the same side with
respect to $\gamma$. Assuming this property and also that
$\gamma$ minimizes the number of intersection points with
$\alpha_1 \cup \dots \cup \alpha_n$ in its isotopy class, we
construct the curves $\gamma_1, \dots, \gamma_n$ as in the proof
of Lemma \ref{betagamma/thm} with the $\alpha_i$'s in place of
the $\beta_i$. In particular we get $\gamma_i = \alpha_i$ if $i
\in S$. At this point, we can carry on the proof analogously to
the proof of Lemma \ref{betagamma/thm}, with the only difference
that, each time a curve $\gamma_i$ with $i \in S$ is involved in
the reasoning, we get a good definition of $\delta$ instead of a
contradiction with respect to the minimality of $\gamma$.
\end{proof}

%{\bf Teorema 3.3.} $L_p$ is generated by a finite set of liftable
%powers of intervals.

\medskip\begin{proof}[Proof of Theorem B]
We proceed by induction on the number $n$ of branch points of $p$.
For $n = 1$ the result is trivial. So, let us suppose $n > 1$.

On the other hand, the case $p \cong p_n$ has been examined in
Lemma \ref{liftablepowers/thm}. Hence we can also assume
$p\not\cong p_n$, in such a way that there exists $m \leq d-1$
minimum index such that $\phi(\alpha_m) = \phi(\alpha_{m+1})$.
Then $p^{\alpha_m}$ is connected and $\phi(\alpha_m) =
\phi(\alpha_{m+1}) = (m\;m{+}1)$.

We start by observing that, if $b \in \L_p$ and there exists a
curve $\alpha$ for $p$ such that $p^\alpha$ is connected and
$(\alpha)b$ is $\I_p$-equivalent to $\alpha$, then $b \in \I_p$.
Indeed, if $c \in \I_p$ is such that $(\alpha)b = (\alpha)c$, then
$(\alpha)bc^{-1} = \alpha$ and therefore $bc^{-1}$ can be thought
as a braid in $\L_{p^\alpha}$. By the induction hypothesis, we
have $bc^{-1} \in \I_{p^\alpha} \subset \I_p$ and therefore $b
\in \I_p$. It is easy to see that an analogous argument also
holds if $p^\alpha$ is not connected but has at most one
non-trivial component.

Now, let $b \in \L_p$ be an arbitrary liftable braid. By Lemma
\ref{betadelta/thm}, the curve $\beta = (\alpha_m)b$ is
$\I_p$-equivalent to a curve $\gamma$ such that $\gamma \cap
\alpha_i = \{*\}$ for some $i = 1, \dots, n$. Moreover, $\gamma$
starts from $*$ on the right of $\alpha_i$ if $i < m$ and on the
left of $\alpha_i$ if $i\geq m$. At this point, we conclude the
proof by checking separately the three possible cases.

(1) $i < m$. In this case $\phi(\alpha_i) = (i\,\,i{+}1)$ and both
the restrictions $p^{\alpha_i}$ and $p^{\alpha_i, \alpha_m}$ have
two components, one of which is trivial (the one corresponding to
the sheet $i + 1$ with respect to the base point $*'$). On the
other hand, $p^\gamma$ is connected and therefore the components
of $p^{\alpha_i, \gamma}$ can not be more than two and they
coincide with the ones of $p^{\alpha_i}$. By Lemma
\ref{systems/thm}, there exists  $c \in \L_p$ such that
$(\alpha_i)c = \alpha_i$ and $(\alpha_m)c = \gamma$. By applying
the induction hypothesis to $c$ thought as a braid in
$\L_{p^{\alpha_i}}$, we have that $c \in \I_{p^{\alpha_i}}
\subset \I_p$ and therefore $\beta = (\alpha_m)b$ is
$\I_p$-equivalent to $\alpha_m$. Finally, the starting
observation enable us to conclude that $b \in \I_p$.

(2) $i = m,m+1$ or $i > m +1$ with $m = d - 1$. In this case the
interval $x \simeq \gamma \cup \alpha_i$ is of type $1$ and
$\gamma$ is $\I_p$-equivalent to $\alpha_i$ and therefore to
$\alpha_m$. Then $b \in \I_p$, since $\beta = (\alpha_m)b$ is
$\I_p$-equivalent to $\alpha_m$.

(3) $i > m+1$ with $m < d-1$. In this case we have
$\phi(\alpha_i) = (l\;l{+}1)$ with $l > m$, moreover the
restrictions $p^{\alpha_i}$ and $p^{\alpha_m, \alpha_i}$ are both
connected or they have two components one of which is trivial
(the one corresponding to the sheet $i+1$ with respect to the
base point $*'$). We consider a fundamental system $\delta_1,
\dots, \delta_{n-2}, \gamma, \alpha_i$ for $p$ and set
$\phi(\delta_j) = \sigma_j$ for each $j = 1, \dots, n-2$. Then
$\sigma_1 \dots \sigma_{n-2} = \phi(\omega)\, (l\;l{+}1)\,
(m\;m{+}1) = (m\;m{-}1\;\dots\;1)\,\sigma (l\;l{+}1)\,(m\;m{+}1)$
with $\sigma$ product of cycles all disjoint from
$(m\;m{-}1\;\dots\;1)$. It follows that $(\sigma_1 \dots
\sigma_{n-2})^m(m) = m+1$. Hence the orbits of the action of
$\langle \sigma_1, \dots, \sigma_{n-2}\rangle \subset \Sigma_d$
coincide with the ones of the action of $\langle \sigma_1, \dots,
\sigma_{n-2}, (m\;m{+}1) \rangle$, so that also the components of
$p^{\gamma, \alpha_i}$ correspond to the ones of $p^{\alpha_i}$.
By Lemma \ref{systems/thm}, there exists $c \in \L_p$ such that
$(\alpha_m)c = \gamma$ and $(\alpha_i)c = \alpha_i$. Then, we can
conclude that $b \in \I_p$ by the same argument of case (1).

At this point, in order to prove that $\L_p$ is finitely generated
and therefore can be generated by a finite set of liftable powers
of intervals, it suffices to observe that $\L_p$ is a subgroup of
finite index of $\B_n$ (see \cite{MKS66}). In fact, given $b,c
\in \B_n$, we have that $bc^{-1}\in \L_p$ if and only if
$\varphi((\alpha_i) bc^{-1})=\varphi(\alpha_i)$ for every
$i=1,\dots,n$, by Lemma \ref{liftability/thm}. Then $b$ and $c$
belong to the same coset of $\L_p$ in $\B_n$ if and only if
$\varphi((\alpha_i)b) = \varphi((\alpha_i)c)$ for every $i = 1,
\dots, n$. This means that there is a bijective correspondence
between cosets of $ \L_p$ in $\B_n$ and admissible sequences of
transpositions of $\Sigma_d$ of length $n$. Therefore $\vert
B_n:L_p\vert\le (d(d-1)/2)^n$.
\end{proof}

\end{document}